\documentclass[11pt]{article}
\usepackage[utf8]{inputenc}
\usepackage[english]{babel}
\usepackage{authblk}
\usepackage{graphicx}        % standard LaTeX graphics tool
\usepackage[caption=false]{subfig}
\usepackage{multicol}        % used for the two-column index
\usepackage[bottom]{footmisc}% places footnotes at page bottom
\usepackage{placeins}
\usepackage{changes,todonotes}
\usepackage{pgf}
\usepackage{psfrag}
\usepackage{pgfplots}
\usepackage{hyperref}
\usepackage{appendix}
\usepackage{amsmath,amssymb}
\usepackage{bm}
\usepackage{mathtools}
\usepackage{color}
\usepackage{stmaryrd}
\usepackage{amsthm}
\usepackage{blindtext}
\usepackage{caption}
\usepackage{multirow}
\usepackage[table]{xcolor}

\usepackage[margin=3cm]{geometry}

\usepackage[ruled,vlined]{algorithm2e}

\graphicspath{./FIGURES/}
\newcommand{\vertiii}[1]{{\left\vert\kern-0.25ex\left\vert\kern-0.25ex\left\vert #1 
    \right\vert\kern-0.25ex\right\vert\kern-0.25ex\right\vert}}

\newcommand{\trinorm}[1]{{\vert\kern-0.25ex\vert\kern-0.25ex\vert #1 \vert\kern-0.25ex\vert\kern-0.25ex\vert}}

\definecolor{farbe}{gray}{0.80}

\newtheorem{remark}{Remark}

\begin{document}

\title{Numerical Verification of PolyDG Algebraic Solvers \\ for the Pseudo-Stress Stokes Problem}

\author[$\star$]{Paola F. Antonietti}
\author[$\star$]{Alessandra Cancrini}
\author[$\star$]{Gabriele Ciaramella}

\affil[$\star$]{MOX, Laboratory for Modeling and Scientific Computing, Dipartimento di Matematica, Politecnico di Milano, Piazza Leonardo da Vinci 32, I-20133 Milano, Italy}
%\affil[$\diamond$]{Correspondig author}

\affil[ ]{\texttt {\{paola.antonietti,alessandra.cancrini,gabriele.ciaramella\}@polimi.it}}

\maketitle

\abstract{This work focuses on the development of efficient solvers for the pseudo-stress formulation of the unsteady Stokes problem, discretised by means of a discontinuous Galerkin method on polytopal grids (PolyDG). The introduction of the pseudo-stress variable is motivated by the growing interest in non-Newtonian flow models and coupled interface problems, where the stress field plays a fundamental role in the physical description. The space-time discretisation of the problem is obtained by combining the PolyDG approach in space with the implicit Euler method for time integration.
The resulting linear system, characterised by a symmetric, positive, definite matrix, exhibits deteriorating convergence with standard solvers as the time step decreases. To address this issue, we investigate two tailored strategies: deflated Conjugate Gradient, which mitigates the effect of the most problematic eigenmodes, and collective Block-Jacobi,  which exploits the block structure of the system matrix. Numerical experiments show that both approaches yield iteration counts effectively independent of $\Delta t$, ensuring robust performance with respect to the time step. Future work will focus on extending this robustness to the spatial discretisation parameter $h$ by integrating multigrid strategies with the time-robust solvers developed in this study.}

%\newpage

\section{Introduction}

The numerical approximation of incompressible viscous flows is a fundamental problem in computational fluid dynamics. In this setting, the pseudo-stress formulation of the unsteady Stokes equations \cite{Cancrini} has attracted attention for its relevance to non-Newtonian models and coupled interface problems, where an accurate representation of the stress tensor is essential.
To discretise the problem, we employ a discontinuous Galerkin method on polytopal meshes (PolyDG) and the implicit Euler method for time integration. At each time step, this discretisation leads to a linear system whose matrix is symmetric and positive definite. When refining the mesh or varying the time step, it is crucial to rely on linear solvers and preconditioners whose performance remains stable with respect to the discretisation parameters.
The main goal of this work is to investigate solver robustness with respect to the time step $\Delta t$. As shown by numerical experiments, both the Conjugate Gradient (CG) method  \cite{Magnus} and the CG preconditioned with standard Block-Jacobi suffer from a severe deterioration of convergence as $\Delta t$ decreases.
This behaviour is due to the conditioning of the system matrix, which scales as $1/\Delta t$; therefore, smaller time steps produce increasingly ill-conditioned systems and, consequently, a larger number of iterations.
%leading to a significant increase in the number of iterations.
The structure of the system matrix motivates the use of two tailored solvers: Deflated CG \cite{Nicolaides, Kahl}, which accelerates convergence by eliminating the most problematic eigenmodes, and collective Block-Jacobi \cite{Ciaramella}, which treats pseudo-stress components collectively at the element level.
Numerical tests show that both approaches yield iteration counts independent of $\Delta t$, thereby achieving the desired robustness. This study represents a first step toward a solver that is robust with respect to both $\Delta t$ and $h$, with future work incorporating multigrid techniques for polytopal discretisations to further enhance robustness under spatial refinement.

\section{Model problem and numerical discretization}\label{sec:stokes_problem}

\subsection{Pseudo-stress weak formulation}

We focus on the Stokes problem, which models incompressible viscous free flows, and formulate it in terms of a pseudo-stress unknown rather than its classical expression. The derivation of the pseudo-stress formulation for the unsteady Stokes problem is presented in \cite{Cancrini}. The pseudo-stress variable is defined as $\boldsymbol{\sigma}(\boldsymbol{u},p) = \mu\boldsymbol{\nabla}\boldsymbol{u} - p\mathbb{I}_d$, where $\boldsymbol{u}$ is the flow velocity and $p$ its pressure. Let $\Omega\subset\mathbb{R}^d$, \textcolor{black}{$d=2$}, be an open, convex \textcolor{black}{polygonal} domain with Lipschitz boundary $\partial\Omega$, and let ${\rm \textbf{dev}}(\boldsymbol{\tau}) = \boldsymbol{\tau} - \frac1d {\rm tr}(\boldsymbol{\tau})\mathbb{I}_d$ be the deviatoric operator, where ${\rm tr}(\cdot)$ is the trace operator. Then, the Stokes problem can be rewritten as

\begin{equation*}
\begin{cases}
  \frac1\mu \frac{\partial {\rm \textbf{dev}}(\boldsymbol{\sigma})}{\partial t}
  - \boldsymbol{\nabla}\left(\boldsymbol{\nabla}\cdot\boldsymbol{\sigma}\right)=\boldsymbol{F}, 
  \qquad &\text{in }\Omega\times(0,T], \\
  \boldsymbol{\nabla}\cdot\boldsymbol{\sigma}=\boldsymbol{g}_D, \;\; &{\rm on}\ \Gamma_D\times(0,T], \\
  \boldsymbol{\sigma}\ \boldsymbol{n} = \boldsymbol{g}_N,\;\; 
  &{\rm on}\ \Gamma_N\times(0,T], \\
  {\rm \textbf{dev}}(\boldsymbol{\sigma})(\cdot,t=0) = \boldsymbol{\sigma}_0, &\text{in }\Omega, 
\end{cases}
\label{eq:eqnn}
\end{equation*}
where $\mu>0$ is the fluid viscosity, and $T>0$ the final simulation time. The boundary of $\Omega$ is partitioned as $ \Gamma_D\cup \Gamma_N = \partial\Omega$, with $\Gamma_D\cap \Gamma_N = \emptyset$. For simplicity, we assume both $\left|  \Gamma_D \right| > 0$ and $\left|   \Gamma_N \right| > 0$. 
%with $|\cdot|$ denoting the Hausdorff measure.
We also assume that $\boldsymbol{F} \in \mathbb{H}^1(\Omega)$, $\boldsymbol{g}_N \in \boldsymbol{L}^2(\Gamma_N)$, $\boldsymbol{\sigma}_0 \in \mathbb{L}^2(\Omega)$ and that $\forall t$ the Dirichlet datum $\boldsymbol{g}_D=\boldsymbol{g}_D(t)$ is the trace of a function in $\boldsymbol{H}^1(\Omega)$. 
To strongly enforce the essential traction condition on $\Gamma_N$, we define 
%the subspace 
$$
\mathbb{H}_{0,\Gamma_N}({\rm div}, \Omega) = \{\boldsymbol{\eta}\in\mathbb{H}({\rm div}, \Omega)\ | \ \langle\bm\eta\ \bm n,\bm v\rangle_{\partial\Omega}=0 \ \ \forall \bm v\in \bm H^1_{0,\Gamma_D}(\Omega)\},
$$
where $\bm H^1_{0,\Gamma_D}(\Omega)= \{\boldsymbol{v} \in \bm H^1(\Omega)^d \;|\; \boldsymbol{v}=\boldsymbol{0} \; \text{on}\; \Gamma_D\}$.
Then, the corresponding weak formulation reads as: for any $t\in (0, T]$, find $\bm \sigma(t) \in \mathbb{H}_{0,\Gamma_N}({\rm div}, \Omega)$ such that  
\begin{equation}\label{weak_stress}
     (\mu^{-1} \partial_t{\rm \textbf{dev}}(\boldsymbol{\sigma}), {\rm \textbf{dev}}(\boldsymbol{\tau}))_{\Omega} + (\boldsymbol{\nabla}\cdot\boldsymbol{\sigma}, \boldsymbol{\nabla}\cdot\boldsymbol{\tau})_{\Omega} =  (\boldsymbol{F}, \boldsymbol{\tau})_{\Omega}  +  \langle\bm g_D, \boldsymbol{\tau}\,\boldsymbol{n}\rangle\textcolor{black}{_{\partial \Omega}}
\end{equation}
for any $\boldsymbol{\tau}\in \mathbb{H}_{0,\Gamma_N}({\rm div}, \Omega)$. For further details, the reader is referred to \cite{Cancrini}.

\subsection{PolyDG semi-discrete formulation}

We can now introduce the PolyDG semi-discrete formulation of \eqref{weak_stress}. Let $\mathcal{T}_h$ be a polytopal mesh of the domain $\Omega$, i.e.,   $\mathcal{T}_h = \textcolor{black}{\bigcup_{\kappa} \kappa}$, being $\kappa$ a  \textcolor{black}{general polygon ($d = 2$).} % or polyhedron ($d = 3$). 
 Given a polytopal element $\kappa$, we define by $\left|  \kappa \right|$ its measure and by \textcolor{black}{$h_\kappa$} its diameter, and set \textcolor{black}{$h = \max_{\kappa \in \mathcal{T}_h} h_\kappa$}. We let a polynomial degree \textcolor{black}{$p_\kappa \ge 1$} be associated with each element $\kappa \in \mathcal{T}_h$ and we denote by $p_h : \mathcal{T}_h \rightarrow \mathbb{N^*} = \{n \in \mathbb{N} : n \ge 1\}$ the piecewise constant function such that $(p_h)|_{\kappa} = \textcolor{black}{p_\kappa}$. Then, we define the discrete space
$\boldsymbol{V}_h = [\textit{P}_{p_h} (\mathcal{T}_h)]^{d\times d}$, where $\textit{P}_{p_h} (\mathcal{T}_h) = \Pi_{\kappa \in \mathcal{T}_h} \mathbb{P}_{p_\kappa} (\kappa)$, and $\mathbb{P}_{\ell} (\kappa)$ is the space of piecewise polynomials in $\kappa$ of total degree less than or equal to $\ell \geq 1.$ 
 We define an interface as the intersection of the ($d - 1$)-dimensional faces of any two neighboring
elements of $\mathcal{T}_h$. 
%If $d = 2$, an interface/face is a line segment and the set of all interfaces/faces is denoted by $\mathcal{F}_h$.
%When $d = 3$, an interface can be a general polygon that we assume could be further decomposed into a set of planar triangles collected in the set $\mathcal{F}_h$.
We also decompose the set of faces as $\mathcal{F} = \mathcal{F}_h^I \cup \mathcal{F}_h^D  \cup \mathcal{F}_h^N$, where $\mathcal{F}_h^I$ contains the internal faces and $\mathcal{F}_h^D$ and $\mathcal{F}_h^N$ the faces of the Dirichlet and Neumann boundary, respectively. We refer the reader to \cite{Cancrini, CangianiDongGeorgoulisHouston_2017} for the main assumptions on $\mathcal{T}_h$. 
Finally, for sufficiently piecewise smooth vector- and tensor-valued fields $\bm{v}$ and $\bm{\tau}$, respectively, and for any pair of neighbouring elements $\kappa^+$ and $\kappa^-$ sharing a face $F\in \mathcal{F}_h^I$,  we introduce the jump and average operators 
$
[[ \bm{v} ]]  = \bm{v}^+\otimes\bm{n}^++\bm{v}^-\otimes\bm{n}^-, [[\boldsymbol{\tau}]] = \boldsymbol{\tau}^+ \boldsymbol{n}^+ +  \boldsymbol{\tau}^- \mathbf{n}^-, 
 \{\!\{\bm{v}\}\!\} = \frac{(\bm{v}^+ + \  \bm{v}^-)}{2},  
 \{\!\{\boldsymbol{\tau}\}\!\} = \frac{(\boldsymbol{\tau}^+  + \  \boldsymbol{\tau}^-)}{2},
 $
where $\otimes$ is the tensor product in $\mathbb{R}^d$, $\cdot^{\pm}$ denotes the trace on $F$ taken within $\kappa^\pm$, and $\bm{n}^\pm$ is the outer normal vector to $\partial \kappa^\pm$. 
Accordingly, on boundary faces  $F \in \mathcal{F}_h^D  \cup \mathcal{F}_h^N$, 
we set $[[\bm{v}]] = \bm{v}\otimes\bm{n}, \,\,
[[\boldsymbol{\tau}]] = \boldsymbol{\tau} \boldsymbol{n}, \,\, {\rm and }  \,\,
\{\!\{\bm{v}\}\!\} =  \bm{v}, \,\,
\{\!\{\bm{\tau}\}\!\} =  \bm{\tau}.$
In the following, we use $\nabla_h \cdot$ to denote the element-wise divergence operator, and we use the short-hand notation $(\cdot,\cdot)_{\mathcal{T}_h} = \sum_{\kappa\in \mathcal{T}_h}\int_{\kappa} \cdot$ and $ \langle \cdot,\cdot \rangle_{\mathcal{F}_h} = \sum_{F\in \mathcal{F}_h}\int_{F} \cdot$.
We consider the following semi-discrete PolyDG approximation to \eqref{weak_stress}: for any $t\in (0,T]$, find $\bm \sigma_h(t) \in \bm V_h$ s.t.
\begin{equation}\label{weak_dg}
   \begin{cases}
\mathcal{M}(\partial_t \bm \sigma_{h}, \bm \tau_h) + \mathcal{A} (\bm \sigma_h, \bm \tau_h) = F(\bm \tau_h) & \forall \, \bm \tau_h \in  \bm V_h,  \\
(\bm \sigma_h(0), \bm \tau_h) = (\bm \sigma_0, \bm \tau_h)  & \forall \, \bm \tau_h \in  \bm V_h,
   \end{cases}
\end{equation}
where for any $\bm \sigma, \bm \tau \in \bm V_h$ we have defined
\begin{eqnarray*}
    \mathcal{M}(\bm \sigma, \bm \tau) & = & ( \mu^{-1} \boldsymbol{{\rm dev}}(\boldsymbol{\sigma}), \boldsymbol{{\rm dev}}(\boldsymbol{\tau}) )_{\mathcal{T}_h}, \\  
    \mathcal{A}(\bm \sigma, \bm \tau) & = & (\boldsymbol{\nabla}_h\cdot\boldsymbol{\sigma}, \boldsymbol{\nabla}_h\cdot\boldsymbol{\tau})_{\mathcal{T}_h} - \langle \{\!\{\boldsymbol{\nabla}_h\cdot\boldsymbol{\sigma}\}\!\} , [[\boldsymbol{\tau}]] \rangle_
    {\mathcal{F}_h^{I,N}} \\
    &  & - \langle \{\!\{\boldsymbol{\nabla}_h\cdot\boldsymbol{\tau}\}\!\} , [[\boldsymbol{\sigma}]] \rangle_
    {\mathcal{F}_h^{I,N}}    
    + \langle \gamma_e[[\boldsymbol{\sigma}]] , [[\boldsymbol{\tau}]] \rangle_
    {\mathcal{F}_h^{I,N}}, \\
    F(\bm \tau) & = &  (\boldsymbol{F}, \boldsymbol{\tau})_{\mathcal{T}_h} +
    \langle \boldsymbol{g}_D,  \boldsymbol{\tau}\boldsymbol{n} \rangle_{\mathcal{F}_h^D}
   +  \langle \boldsymbol{g}_N , \gamma_e \boldsymbol{\tau}\boldsymbol{n} + (\boldsymbol{\nabla}_h\cdot\boldsymbol{\tau})\rangle_{\mathcal{F}_h^N}.
\end{eqnarray*}
%
%We now add "symmetrization" and "stabilization" terms, considering that: 
%
%\begin{itemize}
%\item on $\mathcal{F}_h^I$ we have $$;
%\item on $\mathcal{F}_h^N$ we have $\boldsymbol{\sigma} \cdot \mathbf{n} - \boldsymbol{g}_N = 0$.
%\end{itemize}
%
%Thus, we obtain the following "symmetrization" terms:
%\begin{align*}
%-\theta \sum_{F\in \mathcal{F}_h^I}\int_{F} [[\boldsymbol{\sigma}\boldsymbol{n}]] \cdot  \{\!\{\boldsymbol{\nabla}\cdot\boldsymbol{\tau}\}\!\} \,{\rm d}s &= 0
%\end{align*}
%\begin{align*}
%-\theta \sum_{F\in \mathcal{F}_h^N}\int_{F} \boldsymbol{\sigma}\boldsymbol{n} \cdot (\boldsymbol{\nabla}\cdot\boldsymbol{\tau}) \,{\rm d}s
%+\theta \sum_{F\in \mathcal{F}_h^N}\int_{F} \boldsymbol{g}_N \cdot (\boldsymbol{\nabla}\cdot\boldsymbol{\tau}) \,{\rm d}s &= 0
%\end{align*}
%\\
%and the following "stabilization" terms:
%\begin{align*}
%\gamma \sum_{F\in \mathcal{F}_h^I}\int_{F} [[\boldsymbol{\sigma}\boldsymbol{n}]] \cdot  [[\boldsymbol{\tau}\boldsymbol{n}]] \,{\rm d}s &= 0
%\end{align*}
%\begin{align*}
%\gamma \sum_{F\in \mathcal{F}_h^N}\int_{F} \boldsymbol{\sigma}\boldsymbol{n} \cdot \boldsymbol{\tau}\boldsymbol{n}\,{\rm d}s -
%\gamma \sum_{F\in \mathcal{F}_h^N}\int_{F} \boldsymbol{g}_N \cdot \boldsymbol{\tau}\boldsymbol{n} \,{\rm d}s &= 0
%\end{align*}
Here, $\mathcal{F}_h^{I,N} = \mathcal{F}_h^{I} \cup \mathcal{F}_h^{N}$ and the stabilization function $\gamma_e: \mathcal{F}_h^{I,N}\rightarrow\mathbb{R}_+$ is defined as a function of the penalty coefficient $\alpha>0$ as follows:
\begin{equation*}\label{def:penalty}
    \gamma_e(\bm x) = \begin{cases}
        \alpha \max_{\kappa \in \{ \kappa^+,\kappa^-\}} \frac{p_\kappa^2}{h_\kappa}, & \bm x \in e, e \in \mathcal{F}^I_h, e \subset \partial \kappa^+ \cap \partial \kappa^-,  \\
        \alpha \frac{p_\kappa^2}{h_\kappa}, & \bm x \in e, e \in \mathcal{F}^N_h, e \subset \partial \kappa^+ \cap \partial \Gamma_N.
    \end{cases}
\end{equation*}

\subsection{PolyDG fully-discrete formulation}
 
We introduce a basis $\{ \phi_i, i=1, \dots ,N_h \}$ for the space $\bm V_h$ and express $\boldsymbol{\sigma}_h$ as a linear combination of these basis functions, where the unknown coefficients are functions of time.
We collect the latter in the vector $\bm \sigma_h$, denote by $M$ (resp. $A$) the matrix representation of the bilinear form  $\mathcal{M}(\cdot,\cdot)$ (resp. $\mathcal{A}(\cdot, \cdot)$), and by $\bm f$ the vector representation of the linear functional $F(\cdot)$. The algebraic formulation of \eqref{weak_dg} reduces to: $ \forall t\in(0,T]$, find $\bm \sigma_h(t) \in \bm V_h$ s.t. \ \
$ M\dot{\bm \sigma}_h(t) + A\bm \sigma_h(t) = \bm f(t) \ \  \forall t \in (0,T],$
with initial condition $\bm \sigma_h(0) = \bm \sigma_{0,h}$, being the latter the vector representation of the  $\bm L^2$-projection of $\bm \sigma_0$ onto  $\bm V_h$.  
To integrate this system in time we use the $\theta$-method.
Note that, in this framework, an explicit time-stepping scheme cannot be employed since it would lead to an algebraic system with a singular matrix.
So, we discretize our problem by applying the implicit Euler method (i.e. for $\theta = 1$), and we obtain: for any $n=1,\dots, N_T$ find $\bm \sigma_h^{n+1}$ such that
\begin{equation} \label{pb_study}
A^* \ \bm \sigma_h^{n+1} =  \bm f^*
\end{equation}
with \textcolor{black}{$N_T =T/\Delta t$}, and where the superscript $n$ means the approximation/evaluation of the given quantity at time $t_n = n \Delta t$, $n=0,\dots, N_T$. 
Here, we have introduced the operator $A^* := M + \Delta t A$, together with the corresponding bilinear form $\mathcal{A^*}(\cdot, \cdot) := \mathcal{M}(\cdot, \cdot) + \Delta t \mathcal{A}(\cdot, \cdot)$, which defines a symmetric and positive definite operator. Furthermore, we have defined the right-hand side vector $ \bm f^* := M  \bm \sigma_h^{n} + \Delta t  \bm f^{n+1}$, with the associated operator $ \mathcal{F}^*(\cdot) := \mathcal{M}(\boldsymbol{\sigma}^n,\cdot) + \Delta t \mathcal{F}(\cdot)$.
We also introduce the vector of unknowns
\begin{align} \label{sigma}
{\bm \sigma}=
\begin{bmatrix}
{\bm \sigma}_{11}^\top & 
{\bm \sigma}_{12}^\top & 
{\bm \sigma}_{21}^\top &
{\bm \sigma}_{22}^\top
\end{bmatrix}^\top. 
\end{align}
The matrix $A^*$ in (\ref{pb_study}) is defined in terms of $M$ and $A$, where $M$ is positive semi-definite and $A$ is positive definite. Their structures are given by:
\begin{align*}
M = \begin{bmatrix}
    \frac{1}{2} & 0 & 0 & -\frac{1}{2} \\
    0 & 1 & 0 & 0 \\
    0 & 0 & 1 & 0 \\
    -\frac{1}{2} & 0 & 0 & \frac{1}{2} \\
\end{bmatrix} \otimes M_1, \quad \quad 
A = \begin{bmatrix}
    B_1 & B_2 & 0 & 0 \\
    B_2 & B_3 & 0 & 0 \\
    0 & 0 & B_1 & B_2 \\
    0 & 0 & B_2 & B_3 \\
\end{bmatrix},
\end{align*}
where $M_1, B_1, B_2, B_3
\in \mathbb{R}^{\frac{n}{4}\times\frac{n}{4}}$. In particular, $M_1(i,j) = \int_\Omega \phi_j \phi_i$, while $B_1(i,j) \approx \int_\Omega \partial_x \phi_j\, \partial_x \phi_i$, $B_2(i,j) \approx \int_\Omega \partial_x \phi_j\, \partial_y \phi_i$, and $B_3(i,j)  \approx \int_\Omega \partial_y \phi_j\, \partial_y \phi_i$.

\textcolor{black}{\begin{remark}\label{remark1}
    For completeness, we recall here the main assumptions and results from the error analysis of the fully discrete scheme \cite{Cancrini}. Let $\bm \sigma_h^{n+1} \in \bm V_h$ be the solution of \eqref{pb_study} for a sufficiently large $\alpha$. Then, if the exact solution is sufficiently smooth and $(h, p)$ are quasi uniform, the error in the energy norm satisfies $\| \bm \sigma -\bm \sigma_h \|_{E} \sim \Delta t + h^p$. 
\end{remark}}

\section{Numerical solvers}

In this section, we study numerical solvers for the linear system (\ref{pb_study}) arising from the proposed discretisation. Since the system matrix $A^*$ is symmetric and positive definite, the CG method is the natural solver. However, $A^*$ depends on the singular matrix $M$; therefore, as the time step $\Delta t \to 0$, the conditioning of $A^*$ deteriorates 
and the number of CG iterations increases.
To obtain robustness with respect to $\Delta t$, we consider two strategies: Deflated CG, which removes the influence of the kernel of $M$, and CG preconditioned with collective Block-Jacobi, which exploits the block structure of the matrix to improve convergence. 
%The following sections present these two approaches.

%The theoretical error estimate $\| \bm \sigma -\bm \sigma_h \|_{E}$ and the corresponding numerical validation obtained with this direct approach are reported in~\cite{Cancrini}. Building on these results, we now turn our attention to iterative approaches, exploring Krylov subspace and classical stationary schemes.

\subsection{The Deflated CG method} \label{review_deflated}

Consider a linear system of equations of the form (\ref{pb_study}), 
%\begin{equation} \label{linsys}
%A^* \ \bm x =  \bm f
%\end{equation}
where $A^* \in \mathbb{R}^{n \times n}$ is symmetric and positive definite and $\bm \sigma_h^{n+1}, \ \bm f^* \in \mathbb{R}^{n}$. %In this work we focus on the case in which $A^*$ is a large and sparse matrix; 
To simplify the notation, we denote by $\bm x$ the solution $\bm \sigma_h^{n+1}$ of the linear system. It is well known that the speed of convergence of CG depends on the condition number of the matrix $A^*$ \cite{Magnus, Ciaramella}, and to improve the rate of convergence, it can become mandatory to use a preconditioner. A possible way to precondition is via deflation \cite{Nicolaides, Kahl}. Since the speed of convergence of CG depends on the distribution of the eigenvalues of $A^*$, the idea of deflation is to "hide" parts of the spectrum of $A^*$ from CG such that the CG iteration operates on an equivalent system with a significantly reduced condition number compared to 
the original one. The spectral components that are hidden depend on the chosen deflation subspace $\mathbf{S} \subset  \mathbb{R}^n$, and the resulting improvement in convergence speed is entirely governed by this choice. 
Given $\mathbf{S}^{\perp_{A^*}}=(A^*\mathbf{S})^\perp$ the $A^*-$orthogonal complement of $\mathbf{S}$, it is possible to split the solution $\bm x$ into a component in $\mathbf{S}$ and a component in $\mathbf{S}^{\perp_{A^*}}$ via the $A^*-$orthogonal projection $\pi_{A^*}(\mathbf{S}) \in \mathbb{R}^{n \times n}$ onto $\mathbf{S}$. Given a matrix $V\in \mathbb{R}^{n \times n}$ whose columns form a basis for $\mathbf{S}$, and $\pi_{A^*}(\mathbf{S}) = V(V^\top A^*V)^{-1}V^\top A^*$, we obtain 
\begin{align}\label{x_defl}
   \bm x = (I-\pi_{A^*}(\mathbf{S})) \ \hat{\bm x} + V(V^\top A^*V)^{-1}V^\top  \bm f^*,
\end{align}
where $\hat{\bm x}$ is the solution of the so-called \textit{deflated linear system}
\begin{equation} \label{defl_sys}
    A^*(I-\pi_{A^*}(\mathbf{S})) \ \hat{\bm x} = (I-\pi_{A^*}(\mathbf{S}))^\top \bm f^*.
\end{equation}
Thus,to obtain $\bm x$, one needs to compute the solution $\hat{\bm x}$ of the deflated system. The matrix $A^*(I-\pi_{A^*}(\mathbf{S}))$ is symmetric positive semi-definite \cite{Kahl} and so we can apply the CG method. Matrix singularity poses no obstacle to the standard CG iteration, provided that (\ref{defl_sys}) is consistent, i.e., the right hand side $(I-\pi_{A^*}(\mathbf{S}))^\top\bm f^*$ is in the range of $A^*(I-\pi_{A^*}(\mathbf{S}))$. 
To summarise, one interprets deflated CG as the standard CG algorithm applied to the deflated system (\ref{defl_sys}).
%Therefore, the convergence rate of the deflated CG method can be bounded in terms of the largest and smallest non-zero eigenvalues of $A^*(I-\pi_{A^*}(\mathbf{S}))$.
%We proceed by applying the proposed deflation theory to the algebraic system (\ref{pb_study}).
In order to define the $A^*-$orthogonal projection $\pi_{A^*}(\mathbf{S})$, we introduce $V$ as the kernel of the matrix $M$:
\begin{equation*}
V = \text{ker}(M) = \text{span} \left\{ \frac{1}{\sqrt2}\begin{bmatrix}
    I \\
    0 \\
    0 \\
    I \\ 
\end{bmatrix} \otimes M_1 \right\}.
\end{equation*}
%$
%V = \text{ker}(M) = \text{span} \left\{ \frac{1}{\sqrt2}\begin{bmatrix}
%    I \
%    0 \
%    0 \
%    I \ 
%\end{bmatrix}^\top \otimes M_1 \right\}.
%$ 
This choice is motivated by the fact that the kernel of $M$ captures the problematic spectral components of the system, and deflating this subspace effectively removes the directions responsible for slow convergence of CG. 
Note that applying the deflated CG method requires computing the inverse of $V^\top A^* V$ in (\ref{x_defl}). \textcolor{black}{In the current implementation, the inversion of this operator is performed using a direct solver. 
However, it can be shown that this term is approximately $\tfrac{\Delta t}{2}(B_1 + B_3)$, which 
corresponds to a Laplacian-type operator. As a consequence, multigrid-based solvers could be 
employed to achieve an efficient and scalable inversion. This aspect is beyond the scope of the 
present work, whose primary objective is a preliminary assessment of the proposed approach in 
terms of iteration counts.
}

\subsection{Collective Block-Jacobi as preconditioner for CG}

To solve the linear system (\ref{pb_study}), the standard Block-Jacobi method is commonly used as a preconditioner for the CG method. In this approach, the diagonal blocks are defined by separating each component of the solution vector as in (\ref{sigma}), so that the components are treated independently across the elements. Accordingly, for 
each component $\bm{\sigma}_{ij}$, all elements appear consecutively in the global ordering.
%Each block therefore corresponds to a single component on a single element and can be updated in parallel. 
In contrast, the collective Block-Jacobi approach groups the four components 
associated with a given element into a single block, so that the unknowns belonging to element $\kappa$ 
are kept adjacent. It can be shown that the resulting permuted matrix is block diagonally dominant, in agreement with the definition in \cite{Varga, Price}, and this also explains the improved robustness with respect to $\Delta t$ observed in the numerical tests.

%As we will show in the numerical tests, the standard Block-Jacobi preconditioner is not robust with respect to $\Delta t$, in contrast to the collective Block-Jacobi approach. %This difference is directly related to the block diagonal dominance properties of the corresponding preconditioned matrices.

\section{Numerical results}

In this section, we present numerical experiments aimed at assessing the performance of the proposed strategies for solving the algebraic system (\ref{pb_study}) at each time step, implemented in the open source MATLAB library \texttt{lymph} \cite{lymph2024}. 
%The analysis focuses on evaluating the influence of $\Delta t$ by examining the convergence behavior of various methods while varying both the time step and the mesh size.
 For all numerical tests, the parameters are fixed as follows: the polynomial degree is set to $p=3$, the domain considered is $\Omega = (0,1)^2 $ and the polytopal meshes consist of 50 elements ($h_0 \approx 0.2462$), 100 elements ($h_1 \approx 0.1759$), 200 elements ($h_2 \approx 0.1260$), 400 elements ($h_3 \approx 0.0909$), and 800 elements ($h_4 \approx 0.0637$). \textcolor{black}{The grey-shaded entries in the tables indicate comparable time steps and mesh sizes, and thus, according to Remark~\ref{remark1}, the regime where time and space discretisation errors are balanced.}
First, we compute the condition number of $A^*$ for these five values of $h$ and different time steps $\Delta t$ (left part of Table~\ref{tab:h_dt}), observing that for small $\Delta t$ the condition number is proportional to $1/\Delta t$. The right part of Table~\ref{tab:h_dt} reports the condition number of the matrix preconditioned with the collective Block-Jacobi method, highlighting the improvement in conditioning due to the preconditioner. Additional numerical tests (not shown here) confirm that the condition number of the standard Block–Jacobi preconditioned matrix exhibits the same dependence on $\Delta t$
as that of $A^*$.
%Numerical results confirm that the condition number of the standard Block–Jacobi preconditioned matrix exhibits the same dependence on $\Delta t$ as that of $A^*$.
%
\begin{table}[h!]
\centering
\renewcommand{\arraystretch}{1.2} % aumenta spazio verticale
\setlength{\tabcolsep}{1.6pt}      % aumenta spazio orizzontale tra colonne

\begin{tabular}{c|ccccc|ccccc}
\hline
\textbf{$\Delta t$} & \textbf{$h_0$} & \textbf{$h_1$} & \textbf{$h_2$} & \textbf{$h_3$} & \textbf{$h_4$} & \textbf{$h_0$} & \textbf{$h_1$} & \textbf{$h_2$} & \textbf{$h_3$} & \textbf{$h_4$} \\
\hline
$10^{-2}$ &  \cellcolor{lightgray} 8.5$\cdot 10^{5}$ & 1.6$\cdot 10^{6}$ & 3.0$\cdot 10^{6}$ & 6.5$\cdot 10^{6}$ & 1.3$\cdot 10^{7}$ & \cellcolor{lightgray} 7.2$\cdot 10^{5}$ & 1.2$\cdot 10^{6}$ & 2.7$\cdot 10^{6}$ & 8.0$\cdot 10^{6}$ & 1.3$\cdot 10^{7}$ \\
$10^{-3}$ &  \cellcolor{lightgray} 5.6$\cdot 10^{5}$ & \cellcolor{lightgray} 1.0$\cdot 10^{6}$ & \cellcolor{lightgray} 1.9$\cdot 10^{6}$ & 3.9$\cdot 10^{6}$ & 8.0$\cdot 10^{6}$ & \cellcolor{lightgray} 4.3$\cdot 10^{5}$ & \cellcolor{lightgray} 7.4$\cdot 10^{5}$ & \cellcolor{lightgray} 1.6$\cdot 10^{6}$ & 4.5$\cdot 10^{6}$ & 7.6$\cdot 10^{6}$  \\
$10^{-4}$ &  \cellcolor{lightgray} 5.2$\cdot 10^{5}$ & \cellcolor{lightgray} 9.4$\cdot 10^{5}$ & \cellcolor{lightgray} 1.8$\cdot 10^{6}$ & \cellcolor{lightgray} 3.5$\cdot 10^{6}$ & \cellcolor{lightgray} 7.3$\cdot 10^{6}$ & \cellcolor{lightgray} 3.6$\cdot 10^{5}$ & \cellcolor{lightgray} 6.5$\cdot 10^{5}$ & \cellcolor{lightgray} 1.4$\cdot 10^{6}$ & \cellcolor{lightgray} 4.1$\cdot 10^{6}$ & \cellcolor{lightgray} 6.9$\cdot 10^{6}$  \\
$10^{-5}$ & 5.7$\cdot 10^{5}$ & \cellcolor{lightgray} 9.8$\cdot 10^{5}$ & \cellcolor{lightgray} 1.8$\cdot 10^{6}$ & \cellcolor{lightgray} 3.5$\cdot 10^{6}$ & \cellcolor{lightgray} 7.2$\cdot 10^{6}$ & 2.8$\cdot 10^{5}$ &  \cellcolor{lightgray} 5.2$\cdot 10^{5}$ & \cellcolor{lightgray} 1.1$\cdot 10^{6}$ & \cellcolor{lightgray} 3.4$\cdot 10^{6}$ & \cellcolor{lightgray} 6.3$\cdot 10^{6}$ \\
$10^{-6}$ & 1.7$\cdot 10^{6}$ & 1.9$\cdot 10^{6}$ & 2.5$\cdot 10^{6}$ &  \cellcolor{lightgray} 4.1$\cdot 10^{6}$ & \cellcolor{lightgray} 7.8$\cdot 10^{6}$ & 2.7$\cdot 10^{5}$ & 4.8$\cdot 10^{5}$ & 8.7$\cdot 10^{5}$ & \cellcolor{lightgray} 2.5$\cdot 10^{6}$   & \cellcolor{lightgray} 4.7$\cdot 10^{6}$ \\
$10^{-7}$ & 1.6$\cdot 10^{7}$ & 1.6$\cdot 10^{7}$ & 1.5$\cdot 10^{7}$ & 1.5$\cdot 10^{7}$ & 1.7$\cdot 10^{7}$ & 2.7$\cdot 10^{5}$ & 4.8$\cdot 10^{5}$ & 8.5$\cdot 10^{5}$ & 2.3$\cdot 10^{6}$ & 4.1$\cdot 10^{6}$ \\
$10^{-8}$ & 1.6$\cdot 10^{8}$ & 1.6$\cdot 10^{8}$ & 1.5$\cdot 10^{8}$ & 1.4$\cdot 10^{8}$ & 1.5$\cdot 10^{8}$ & 2.7$\cdot 10^{5}$ & 4.7$\cdot 10^{5}$ & 8.5$\cdot 10^{5}$ & 2.2$\cdot 10^{6}$ & 4.1$\cdot 10^{6}$ \\
$10^{-9}$ & 1.6$\cdot 10^{9}$ & 1.6$\cdot 10^{9}$ & 1.5$\cdot 10^{9}$ & 1.4$\cdot 10^{9}$ & 1.4$\cdot 10^{9}$ & 2.7$\cdot 10^{5}$ & 4.7$\cdot 10^{5}$ & 8.5$\cdot 10^{5}$ & 2.2$\cdot 10^{6}$ & 4.0$\cdot 10^{6}$  \\
$10^{-10}$ & 1.6$\cdot 10^{10}$ & 1.6$\cdot 10^{10}$ & 1.5$\cdot 10^{10}$ & 1.4$\cdot 10^{10}$ & 1.4$\cdot 10^{10}$ & 2.7$\cdot 10^{5}$ & 4.7$\cdot 10^{5}$ & 8.5$\cdot 10^{5}$ & 2.2$\cdot 10^{6}$ & 4.0$\cdot 10^{6}$ \\
\hline
\end{tabular}

\caption{Condition number of the matrix $A^*$ (left) and condition number of the matrix preconditioned with collective Block-Jacobi (right), as functions of the mesh size $h$ and the time step $\Delta t$.}
\label{tab:h_dt}
\end{table}

\noindent Table~\ref{tab:DCGsolvers} and Table~\ref{tab:BJsolvers} compare the convergence behavior of the standard CG, Deflated CG, and CG preconditioned with standard and collective Block-Jacobi for different time steps 
$\Delta t$ and mesh sizes $h$. For each method, the tables report the number of iterations required to reduce the (preconditioned) relative residual below $10^{-8}$, considering only the 
 solution of the linear system (\ref{pb_study}) arising at a single time step. The iteration counts are obtained as the average over 10 independent tests, in which the initial condition is perturbed using randomly generated data. In particular, Table~\ref{tab:DCGsolvers} shows that the standard CG iteration counts grow rapidly for small $\Delta t$ and refined meshes, while the Deflated CG method yields dramatically fewer iterations, confirming the effectiveness of the deflation strategy in improving convergence.
 
 %In particular, Table~\ref{tab:DCGsolvers} shows that standard CG iterations increase sharply for small $\Delta t$ and fine meshes, while Deflated CG significantly reduces iteration counts, demonstrating the effectiveness of deflation in improving convergence.

\begin{table}[h!]
\centering
\renewcommand{\arraystretch}{1}
\setlength{\tabcolsep}{6pt}

% --- Prima riga di tabelle (CG e Deflated CG) ---
%\begin{minipage}{0.4\textwidth}
\centering
\begin{tabular}{c|cccc|cccc}
\hline
$\Delta t$ & $h_1$ & $h_2$ & $h_3$ & $h_4$ & $h_1$ & $h_2$ & $h_3$ & $h_4$ \\
\hline
$10^{-2}$    & 1702& 2442 & 3439 & 4727 & 1170 & 1627 & 2256 & 3188 \\
$10^{-3}$   & \cellcolor{lightgray} 1348 & \cellcolor{lightgray} 1949 & 2542 & 3559 & \cellcolor{lightgray} 492 & \cellcolor{lightgray} 660 & 912 & 1286 \\
$10^{-4}$  & \cellcolor{lightgray} 1291 & \cellcolor{lightgray} 1816 & \cellcolor{lightgray} 2329 & \cellcolor{lightgray} 3264 & \cellcolor{lightgray} 221 & \cellcolor{lightgray} 285 & \cellcolor{lightgray} 390 & \cellcolor{lightgray} 525 \\
$10^{-5}$ & \cellcolor{lightgray} 1344 & \cellcolor{lightgray} 1848 & \cellcolor{lightgray} 2339 & \cellcolor{lightgray} 3236 & \cellcolor{lightgray} 106 & \cellcolor{lightgray} 133 & \cellcolor{lightgray} 181 & \cellcolor{lightgray} 220 \\
$10^{-6}$ & 2011 & 2322 & \cellcolor{lightgray} 2640 & \cellcolor{lightgray} 3421 & 62 & 65 & \cellcolor{lightgray} 85 & \cellcolor{lightgray} 100 \\
$10^{-7}$ & 5020 & 5228 & 5367 & 5356 & 61 & 56 & 62 & 57  \\
$10^{-8}$ & 9020 & 10993 & 13665 & 14888 & 60 & 56 & 62 & 58  \\
\hline
\end{tabular}

\caption{Iteration counts for Conjugate Gradient (left) and Deflated CG (right) as functions of $h$ and $\Delta t$, with $\text{tol}=10^{-8}$.}
\label{tab:DCGsolvers}
\end{table}

\noindent In Table~\ref{tab:BJsolvers} we observe that CG preconditioned with standard Block-Jacobi considerably reduces iteration counts compared to the unpreconditioned CG, but the number of iterations still depends on both $\Delta t$ and $h$. In contrast, CG preconditioned with collective Block-Jacobi exhibits nearly uniform iteration counts, independent of $\Delta t$.

\begin{table}[h!]
\centering
\renewcommand{\arraystretch}{1}
\setlength{\tabcolsep}{6pt}
% --- Seconda riga di tabelle (BJ e CS) ---
%\begin{minipage}{0.4\textwidth}
\centering
\begin{tabular}{c|cccc|cccc}
\hline
$\Delta t$ & $h_1$ & $h_2$ & $h_3$ & $h_4$ & $h_1$ & $h_2$ & $h_3$ & $h_4$ \\
\hline
$10^{-2}$    & 848 & 1152 & 1653 & 2207 & 630 & 871 & 1242 & 1680 \\
$10^{-3}$   & \cellcolor{lightgray} 665 & \cellcolor{lightgray} 917 & 1272 & 1624 & \cellcolor{lightgray} 498 & \cellcolor{lightgray} 673 & 959 & 1205 \\
$10^{-4}$  & \cellcolor{lightgray} 637 & \cellcolor{lightgray} 860 & \cellcolor{lightgray} 1081 & \cellcolor{lightgray} 1494 & \cellcolor{lightgray} 465 & \cellcolor{lightgray} 635 & \cellcolor{lightgray} 836 & \cellcolor{lightgray} 1108 \\
$10^{-5}$ & \cellcolor{lightgray} 784 & \cellcolor{lightgray} 987 & \cellcolor{lightgray} 1149 & \cellcolor{lightgray} 1521 & \cellcolor{lightgray} 448 & \cellcolor{lightgray} 615 & \cellcolor{lightgray} 768 & \cellcolor{lightgray} 1075 \\
$10^{-6}$ & 1443 & 1605 & \cellcolor{lightgray} 1799 & \cellcolor{lightgray} 2058 & 441 & 613 & \cellcolor{lightgray} 762 & \cellcolor{lightgray} 1056 \\
$10^{-7}$ & 1960 & 2670 & 3398 & 4071 & 443 & 617 & 772 & 1067 \\
$10^{-8}$ & 1745 & 2466 & 3524 & 4835 & 444 & 616 & 782 & 1074\\
\hline
\end{tabular}

\caption{Iteration counts for CG preconditioned with standard Block-Jacobi (left) and CG preconditioned with collective Block-Jacobi (right) as functions of $h$ and $\Delta t$, with $\text{tol}=10^{-8}$.}
\label{tab:BJsolvers}

\end{table}

\noindent Both Deflated CG and CG preconditioned with collective Block-Jacobi are robust with respect to $\Delta t$, although Deflated CG requires inverting a Laplacian-type operator during setup. To achieve robustness with respect to the mesh size $h$ as well, the natural next step is to combine and adapt these strategies within a multigrid framework.

%Iteration counts for both methods still increase with mesh refinement, indicating that combining these strategies within a multigrid framework could further improve robustness with respect to $h$.

%While the collective Block-Jacobi achieves robustness with respect to $\Delta t$, an increase in iterations is still observed as the mesh is refined. To achieve robustness also with respect to the mesh size $h$, the natural next step is to combine and adapt this strategy within a multigrid framework.

\section*{Acknowledgements}
This work is funded
by the European Union (ERC SyG, NEMESIS, project number 101115663). Views and opinions expressed are however
those of the authors only and do not necessarily reflect those of the European Union or the European Research Council Executive Agency. 
All authors are members of the Indam GNCS group.

%\input{references}
%\bibliographystyle{spmpsci}
%\bibliography{references}

\end{document}